\documentclass[12pt]{article}
    \usepackage{amsmath}
    \usepackage{amssymb}
    \usepackage{amsthm}

\usepackage{pstricks}

\usepackage[T2A]{fontenc}
\usepackage[cp1251]{inputenc}

\usepackage{graphicx,psfrag}

      \newcommand {\al}   {\alpha}          \newcommand {\bt}  {\beta}
                \newcommand {\Gam} {\Gamma}

                 \newcommand {\vphi} {\varphi}

      \newcommand {\pl}   {\partial}        
           
      \newcommand {\RRR}  {{\mathbb R}}     \newcommand {\SSS}  {\mathcal{S}}

        \newcommand {\NNN}  {{\mathbb{N}}}  
              
      \newcommand {\FFF}  {{\cal F}}        
       \newcommand {\IIII}  {{\cal I}}      \newcommand {\rad}  {r}
       \newcommand {\phio}  {\phi}         \newcommand {\kappao}  {k}
       \newcommand {\fff}  {\mathfrak{F}}   \newcommand {\kappat}  {\tau}

     \newcommand {\interval}{[-\pi/2,\, \pi/2]}
     \newcommand {\beq}  {\begin{equation}}
      \newcommand {\eeq}  {\end{equation}}
     \newcommand {\beqo}  {\begin{equation*}}
      \newcommand {\eeqo}  {\end{equation*}}

     \newcommand {\RRRR}  {\mathcal{R}}
     \newcommand {\Sigmapm}  {(S^{d-1} \times S_\rad^{d-1})_\pm}
     \newcommand {\Sigmap}  {(S^{d-1} \times S_\rad^{d-1})_+}
     \newcommand {\Sigmam}  {(S^{d-1} \times S_\rad^{d-1})_-}
     \newcommand {\SigmaPM}  {(S^{d-1} \times S_R^{d-1})_\pm}
     \newcommand {\SigmaP}  {(S^{d-1} \times S_R^{d-1})_+}
     \newcommand {\SigmaM}  {(S^{d-1} \times S_R^{d-1})_-}

      \newtheorem{theorem}{Theorem}

      \newtheorem{zam}{Remark}

\title{The problem of camouflaging via mirror reflections}

 \author{Alexander Plakhov\thanks{Center for R\&{}D in Mathematics and Applications, Department of Mathematics, University of Aveiro, Portugal and Institute for Information Transmission Problems, Moscow, Russia.}}

\begin{document}

\maketitle

\begin{abstract}
This work is related to billiards and their applications in geometric optics. It is known that perfectly invisible bodies with mirror surface do not exist. It is natural to search for bodies that are, in a sense, close to invisible. We introduce a {\it visibility index} of a body measuring the mean angle of deviation of incident light rays, and derive a lower estimate to this index. This estimate is a function of the body's volume and of the minimal radius of a ball containing the body. This result is far from being final and opens a possibility for further research.
\end{abstract}

\begin{quote}
\end{quote}

\begin{quote}
\end{quote}

\section{Introduction}

The idea of invisibility has always been attractive for the people. Stories on magic cap and cloak of invisibility form an essential part of folklore, myths and fairy tales. Various methods of camouflaging establishments and objects of importance have been developed by military in all times; one of the most famous developments of the 20th century is the Stealth technology aiming at making airplanes invisible for radars of the enemy.

In the last decades intensive work is being carried on developing technology of meta-materials possessing unusual properties (see, e.g., \cite{pendry}) having in mind, in particular, creating something like a metamaterial cover which is transparent and possesses unusual refracting properties which make invisible every object placed inside.

An important and interesting mathematical construction in the 2D case is proposed in the paper by Leonhardt \cite{leon}. It describes how to make an object invisible by wrapping a lens (a transparent material with varying refraction index) around it.

There is an interesting question, to what extent can one create the effect of invisibility, if only mirror systems are allowed to use. Mirrors (even curved mirrors) are much easier and cheaper for fabrication than hypothetical meta-material devices, and even than lenses with varying refraction index. Some results in this direction have already been obtained. There exist and are described (connected) bodies invisible from 1 point \cite{invis1point,ebook} and (infinitely connected) invisible from 2 points \cite{invis2points}. There exist (connected and even simply connected) bodies invisible in 1 direction (that is, from an infinitely distant point) \cite{0-resist},  (finitely connected) bodies invisible in 2 directions \cite{invisibility}, as well as (infinitely connected) bodies invisible in 3 \cite{invis3dir} and (in the 2D case) in $n$ directions, where the number $n \in \NNN$ of directions is arbitrary \cite{n_points}.

On the other hand, there are negative results revealing restricted possibilities of mirror systems as compared with more sophisticated technologies. In particular, non-existence of perfectly invisible bodies (that is those that are invisible in any direction or (equivalently) from any point outside the body) is proved. Moreover, a conjecture \cite{ebook} has recently been proposed, stating that the set of light rays that are invisible for any fixed body has measure zero. This conjecture is closely connected with the (similar) Ivrii conjecture \cite{ivrii} on the measure of the set of periodic billiard trajectories in a bounded domain. If Ivrii's conjecture is true then, most probably, true also is the conjecture on invisible light rays.

In real life quite common is the situation when perfect invisibility is impossible to achieve. In such cases one tries to reach the effect of partial invisibility, or camouflaging, when the object, though not disappearing completely, still becomes difficult to detect by an observer. It is natural to set such a question in our problem of mirror invisibility. In order to state the mathematical problem, one needs first to elaborate an index of visibility, a certain positive quantity which is close to zero if the body is, in a natural sense, difficult to detect. This quantity should never vanish, since perfectly invisible bodies do not exist.

Then one needs to consider the question, how small can this index be made in a certain class of bodies. For example, if even the index does not vanish, is it possible to construct a sequence of bodies of constant volume with the index going to zero?

Choosing the visibility index is a difficult task; it is more difficult than just defining the notion of invisibility. The body is observed against a certain background, and the choice will depend, in particular, on the distance of the body from this background. In the limit, when the background is infinitely distant, the visibility index is determined by the angle of deviation of light rays from their original direction and does not depend on transverse displacement of the rays. This limit will be used later on in this paper.

The aim of the paper is to give (partial) answers to the questions stated above. If the body has the volume $A$ and is contained in a sphere of radius  $R$, then its visibility index is not less than a certain positive value, a function of $A$ and $R$. This function goes to zero when $A$ is constant as $R \to \infty$.

\section{Main definitions and statement of the results}

First of all fix the notation. A body with specular surface is a bounded finitely connected domain with piecewise smooth boundary in Euclidean space $\RRR^d$, with $d \ge 2$. It will be called a {\it domain} and designated by $D$. Since everything is about specular reflections in the framework of geometric optics, we adopt the notation of billiard theory and consider motion of billiard particles outside $D$, rather than light rays.

Fix a domain $D$ and take a sphere $S^{d-1}_R$ of radius $R > 0$ centered at the origin and containing $D$. It is assumed that the background lies on the sphere. As a result of observation of the background one must conclude whether the body is or is not present here. For a point $\xi$ on the sphere and a unit vector $v$ such that $\langle v, \xi \rangle < 0$ consider the trajectory of a billiard particle starting at $\xi$ with the velocity $v$ and the half-line with the endpoint $\xi$ and the directing vector $v$, and denote by $\theta = \theta_R(v,\xi)$ the angular distance between the (second) points of intersection of the trajectory and of the half-line with the sphere. We define the measure spaces
$$
\SigmaPM = \{ (v,\xi) \in S^{d-1} \times S^{d-1}_R : \pm\langle v, \xi \rangle \ge 0 \}
$$
with the measure $\mu_R = \mu$ defined by $d\mu(v,\xi) = |\langle v,\, n(\xi) \rangle|\, dv\, d\xi$, where $n(\xi) = \xi/R$ is the outer unit normal to the sphere $S_R^{d-1}$ at the point $\xi$. Take a monotone increasing function $f : [0,\, \pi] \to \RRR$ such that $f(0) = 0$ and consider the value
$$
\FFF_R(D) = \int_{\SigmaM} f(\theta_R(v,\xi))\, d\mu(v,\xi),
$$
which will be called the {\it visibility index} of $D$. Note that the spaces $\SigmaM$ and $\SigmaP$ correspond to billiard trajectories coming in the sphere $S_R^{d-1}$ and going out of it, respectively, and $\mu$ is a natural measure counting the amount of billiard trajectories intersecting the sphere.

We stress that the equality $\FFF_R(D) = 0$ does not yet guarantee invisibility of $D$. In fact the domain is invisible, if and only if $\FFF_R(D) = 0$ for any $R$ sufficiently large.

In the limit $R \to \infty$ the quantity $\theta$ does not depend on the transverse displacement (shift) of the trajectory going away, but only on its direction $v^+$. Let us introduce some more notation. For the billiard trajectory entering the sphere $S_R^{d-1}$ at a point $\xi$ and having a velocity $v$ at this point, we denote by
$$
v^+ = v^+_{D}(v,\xi), \quad \xi^+ = \xi^+_{D,R}(v,\xi)
$$
the second point of intersection of this trajectory with the sphere and its velocity at this point. (Note that the velocity $v^+$ does not depend on the radius $R$ of the sphere containing $D$.) The mapping
$$
T = T_{D,R} : (v, \xi) \mapsto \big(v^+_{D}(v,\xi), \xi^+_{D,R}(v,\xi)\big)
$$
is a one-to-one mapping from $\SigmaM$ onto $\SigmaP$ (defined up to a subset of zero measure) and preserves the measure $\mu$.

Then in the limit mentioned above, $\theta$ is the angle between $v$ and $v^+$,\, $\theta = \arccos\langle v, v^+\rangle$, and one comes to the following formula for the visibility index, $\FFF(D) = \lim_{R\to\infty} \FFF_R(D)$:
$$
\FFF(D) = \int_{\Sigmam} f\big(\arccos \langle v, v^+_D(v,\xi) \rangle\big)\, d\mu(v,\xi),
$$
with $\rad$ being taken sufficiently large. The value of the integral in this formula does not depend on $\rad$.
The visibility index in this case is related to the situation when the distance to the background is much greater than the size of the domain itself.

Denote by $s_{d-1} = |S^{d-1}| = \frac{2\pi^{d/2}}{\Gam(d/2)}$ the area of the $(d-1)$-dimensional unit sphere, and by $b_{d} = \frac{2\pi^{d/2}}{d\Gam(d/2)}$ the volume of the $d$-dimensional ball. One has, in particular, $s_0 = 2,\; s_1 = 2\pi,\; s_2 = 4\pi,\; b_1 = 2,\; b_2 = \pi,\; b_3 = 4\pi/3$.

Assume that
\beq\label{f}
f(\phio) = c\,\phio^\kappao (1 + o(1)) \quad \text{as} \ \, \phio \to 0^+ \quad (\text{with} \ \, c>0,\, \kappao > 0)
\eeq
and introduce the notation
$$
c_d = \left\{
\begin{array}{ll}
c\, \frac{\kappao^\kappao}{(\kappao+1)^{\kappao+1}}\ \frac{\pi}{2^\kappao}, & \text{if}\, \ d = 2\\
c\, \frac{\kappao^\kappao}{(\kappao+1)^{\kappao+1}}\ \frac{s_{d-1}^{\kappao+1} 2^{-1-d\kappao+2\kappao}}{\big( b_{d-1} s_{d-2} B\big( \frac{d-1}{2}, \frac{d-2}{2} \big) \big)^\kappao}, & \text{if}\, \ d \ge 3.
\end{array}
\right.
$$
In particular, in the 3D case we have $c_3 = c\, \frac{\kappao^\kappao}{(\kappao+1)^{\kappao+1}}\, \frac{1}{(2\pi)^{\kappao-1}}$.

The following Theorem \ref{t_main} establishes a connection between the visibility index of a domain, its volume, and the radius of a ball containing this domain.

\begin{theorem}\label{t_main}
Let a domain $D \subset \RRR^d$ be contained in a ball of radius $\rad$. Then its visibility index $\FFF(D)$, its volume $|D|$, and $\rad$ are related by the inequality
$$
\frac{\FFF(D)}{\rad^{d-1}} \ge h_d( |D|/\rad^d ),
 $$
where $h_d$ are functions of a positive variable satisfying
$$
h_d(x) = c_d\, x^{\kappao+1}\, (1 + o(1))\, \ \text{as} \ \, x \to 0^+.
$$
\end{theorem}

Note that the values $\FFF(D)/\rad^{d-1}$ and $|D|/\rad^d$ are preserved under a scaling transformation (applied to both $D$ and the ambient ball). It is natural therefore that Theorem \ref{t_main} relates these two values.

\begin{zam}\label{z1}
It follows from Theorem \ref{t_main} that
$$
\inf_{|D| = \text{const}} \FFF(D) \ge c_d\, |D|^{\kappao+1} \frac{1}{\rad^{\kappao d + 1}}(1 + o(1)) \quad \text{as}\, \ \rad \to \infty.
$$
This means that the infimum of visibility index in the class of domains with fixed volume contained in a certain ball is greater than a positive constant. This constant goes to zero when the radius of the ball tends to infinity. It would be interesting to learn something about the upper bound for this infimum, at least for a special kind of the function $f$ defining the visibility index. In particular, does the infimum go to zero as $\rad \to \infty$? In other words, is it possible to construct a sequence of domains with fixed volume and with the visibility index going to zero?\footnote{In this case the diameter of the bodies should go to infinity.} We do not know the answer to this question.
\end{zam}

One may wish to have a more direct estimate of the visibility index (without the term $o(1)$). We shall derive such estimates in the 2D and 3D cases and for a particular choice of the function $f$. Namely, take $f(\theta) = 1 - \cos \theta$; the resulting visibility index
$$
\fff(D) = \int_{\Sigmam} \big(1 - \langle v, v^+_D(v,\xi) \rangle\big)\, d\mu(v,\xi)
$$
(with $\rad$ taken sufficiently large) has a simple mechanical interpretation: it is just the mean value over all $v$ of the aerodynamic resistance of $D$ in the direction $v$.\footnote{Of course it does not depend on the radius $\rad$ of the ambient sphere and does not change when the domain is displaced inside this sphere. This is proved in the book \cite{ebook} in Proposition 1.1 of Chapter 1.} We shall briefly call it the {\it mean resistance} of $D$.

Note that the mean resistance of a {\it convex} domain $C$ can easily be determined; denoting by $|\pl C|$ the ($(d-1)$-dimensional) area of its boundary, one has $\fff(C) = \frac{4}{d+1}\, b_{d-1} |\pl C|$. In the 2D and 3D cases one has, respectively, $\fff(C) = 8 |\pl C|/3$ and $\fff(C) = \pi |\pl C|$.  In particular, the mean resistances of the 2D and 3D balls, $B_r^2$ and $B_r^3$, of radius $r$ are equal, respectively, to $\fff(B_r^2) = \frac{16}{3}\, \pi r$ and $\fff(B_r^3) = 4\pi^2 r^2$.

Notice the following formulae for the mean resistance in the 2D and 3D cases that are obtained from Theorem \ref{t_main} by direct substitution $c = 1/2,\, \kappao = 2$,
$$
\frac{\fff(D)}{\rad} \ge h_2( {|D|}/{\rad^2} ),\quad \text{where} \ h_2(x) = \frac{\pi}{54}\, x^3 (1 + o(1)),\, \ x \to 0^+ \quad \text{for}\, \ d=2;
$$
$$
\frac{\fff(D)}{\rad^2} \ge h_3( {|D|}/{\rad^3} ),\quad \text{where} \ h_3(x) = \frac{1}{27\pi}\, x^3 (1 + o(1)),\, \ x \to 0^+ \quad \text{for}\, \ d=3.
$$

The following theorem allows one to get rid of the term $o(1)$ in the above formulae.

\begin{theorem}\label{t2}
(a) Let a planar domain $D$ with the area $|D|$ be contained in a circle of radius $\rad$. Then 
$$
\frac{\fff(D)}{\rad} \ge \frac{\pi}{54}\, \bigg( \frac{|D|}{\rad^2} \bigg)^3.
$$
(b) Let a 3-dimensional domain $D$ with the volume $|D|$ be contained in a ball of radius $\rad$. Then 
$$
\frac{\fff(D)}{\rad^2} \ge \frac{1}{27\pi}\, \bigg( \frac{|D|}{\rad^3} \bigg)^3.
$$
\end{theorem}

It is instructive to rewrite these formulas in terms of reduced volume $\kappa_D$ and reduced resistance $\hat\fff_D$ defined by
$$
\kappa_D = \frac{|D|}{b_d \rad^d}, \quad \quad \hat\fff_D = \frac{\fff(D)}{\frac{4}{d+1}\, b_{d-1} s_{d-1} \rad^{d-1}},
$$
where $\rad$ is the radius of the smallest ball containing $D$. One always has $0 < \kappa_D \le 1$, and $\kappa_D = 1$ iff $D$ is a ball. In the latter case we have $\hat\fff_D = 1$. In the case $d=2$ one has $\kappa_D = |D|/(\pi\rad^2)$ and $\hat\fff_D = \fff(D)/(\frac{16}{3}\, \pi\rad)$, and in the case $d=3$ one has $\kappa_D = |D|/(\frac 43\, \pi\rad^3)$ and $\hat\fff_D = \fff(D)/(4\pi^2 \rad^2)$. One can always ensure that the volume of each domain in the sequence is equal to a certain value by taking off the domain a smaller concentric ball.

It is interesting to note that
$$
\sup_{\kappa_D = \kappa} \hat\fff_D = \frac{d+1}{2}
$$
for all $0 < \kappa < 1$. This can easily be derived from Theorem 6.2 in Chapter 6 of the book \cite{ebook} by taking a sequence of domains inscribed in a certain ball and having the property of almost retro-reflection (that is, the initial velocity of the most part of incident particles is reversed as a result of reflections).

On the contrary, only estimates are known for the infimum of $\hat\fff_D$. In particular, the following estimates in the 2D and 3D cases follow directly from Theorem \ref{t2},
\beq\label{for1}
\hat\fff_D \ge \frac{\pi^3}{288}\ \kappa_D^3 \ \ \ \text{\rm for}\, \ d = 2 \quad \quad\ \text{\rm and} \quad \quad\
\hat\fff_D \ge \frac{16}{729}\ \kappa_D^3 \ \ \ \text{\rm for}\, \ d = 3.
\eeq

These estimates are far from being sharp. Indeed, from the same Theorem 6.2 in \cite{ebook} one can derive the exact value of the lower limit of $\hat\fff_D$ when $\kappa_D \to 1$, in particular,
$$
\lim_{\kappat\to 1^-} \inf_{\kappa_{D} = \kappat} \hat\fff_D = m_2 \approx 0.987820 \ \ \text{for}\, \ d=2 \ \text{ and } \
\lim_{\kappat\to 1^-} \inf_{\kappa_{D} = \kappat} \hat\fff_D = m_3 \approx 0.969445 \ \ \text{for}\, \ d=3.
$$
On the other hand, the values of the lower limits given by formulae (\ref{for1}) are much smaller, $\lim_{\kappat\to 1^-} \inf_{\kappa_{D} = \kappat} \hat\fff_D \ge \pi^3/288 \approx 0.11$ for $d=2$ and $\lim_{\kappat\to 1^-} \inf_{\kappa_D = \kappat} \hat\fff_D \ge 16/729 \approx 0.022$ for $d=3$.

There is a question on a natural generalization of formula (\ref{for1}). Consider the relative volume of a domain in its convex envelope, and let the normalized resistance be chosen so that the resistance of the convex envelope of the body equals 1. Is it possible to derive a sensible estimate for the normalized resistance in the spirit of formula (\ref{for1})?

Finally note that the statements of Theorems \ref{t_main} and \ref{t2} hold for broader classes of dynamical systems than billiards. It suffices that the system satisfies the following conditions:
\vspace{1mm}

(i) the motion is free outside a sphere of radius $\rad$;
\vspace{1mm}

(ii) all the trajectories of the system are continuous curves;
\vspace{1mm}

(iii) the natural Liouville measure $dv\, dx$ is invariant under the dynamics of the system.
\vspace{1mm}

One can, for example, take the dynamics with a pseudo-billiard law of reflection off the boundary of $D \subset \RRR^2$; this law is induced by any one-to-one mapping of the segment $\interval$ onto itself preserving the measure $\cos\vphi\, d\vphi$.

\section{Proofs of the theorems}

All statements below are true up to subsets of zero measure.


Fix a value $0 < \phio < \pi$ and let $\Sigma_\phio$ be the set of values $(v,\xi) \in \Sigmam$ such that the angle between $v$ and $v^+ = v^+_D(v,\xi)$ is greater or equal than $\phio$; that is,
$$
\Sigma_\phio = \{ (v,\xi) \in \Sigmam : \langle v, v^+ \rangle \le \cos\phio \}.
$$

We consider the billiard inside the ball $B_\rad^d$ of radius $\rad$ and outside $D$. A particle starts moving at a point $\xi \in S_\rad^{d-1}$ and with a velocity $v \in S^{d-1}$ directed inside the sphere $S^{d-1}$, then makes several reflections off $D$, and finally intersects $S^{d-1}$ again (at the point $\xi^+$ and with the velocity $v^+$) and disappears at the moment of intersection. The phase space of the billiard is $S^{d-1} \times (B_\rad^d \setminus D)$, and its volume $V$ (with respect to the standard Liouville measure $dv\, dx$) equals
\beq\label{eq2}
V = |S^{d-1} \times (B_\rad^d \setminus D)| = s_{d-1} (b_d \rad^d - |D|).
\eeq

On the other hand, denoting by $l(v,\xi)$ the length of the billiard trajectory with the initial data $(v,\xi) \in \Sigmam$ until the final intersection with $S_\rad^{d-1}$ and writing $\xi^+$ in place of $\xi_{D,\rad}^+(v,\xi)$ for brevity, using Santal\'o-Stoyanov formula (see, e.g., \cite{stoyanov}) and the obvious inequality $l(v,\xi) \ge |\xi - \xi^+|$, we have
\beq\label{Vge}
V \ge \int_{\Sigmam} l(v,\xi)\, d\mu(v,\xi) \ge \int_{\Sigmam} |\xi - \xi^+|\, d\mu(v,\xi);
\eeq
(note that a part of the phase space may be inaccessible for particles starting at the boundary of $\Sigmam$).

Further, take an orthonormal coordinate system, $x = (x_1,\ldots x_d)$, centered at the origin and denote by $\RRRR_v$ the rotation of $\RRR^d$ about the origin such that

(a) under this rotation, $v$ goes to $(\bar 0,\, 1) := (0,\ldots,0, 1)$; that is, $\RRRR_v v = (\bar 0, 1)$;

(b) $\RRRR_{v}$ acts as rotation on the 2-dimensional subspace spanned by the vectors $v$ and $(\bar 0, 1)$;

(c) $\RRRR_{v}$ acts as identity on the orthogonal complement to this subspace.

Denote the upper and lower hemispheres of radius $\rad$ by
$$
\SSS_\rad^\pm = \{ \eta = (\eta_1,\ldots, \eta_d) \in S_\rad^{d-1} : \pm \eta_d \ge 0 \}.
$$
Denote $\eta' = (\eta_1,\ldots, \eta_{d-1})$ and consider the standard measure $d\eta'$ on both the hemispheres (that is, the measure of a Borel subset in $\SSS_\rad^\pm$ is the Lebesgue measure of its orthogonal projection on the subspace $\eta_d = 0$).
For each choice of the sign $"+"$ or $"-"$ and for all $v \in S^{d-1}$, the mapping $\xi \mapsto \RRRR_v \xi$ from the hemisphere $\{ \xi \in S^{d-1}_\rad : \pm \langle v, \xi \rangle \ge 0 \}$ with the measure $|\langle v, \xi \rangle|\, d\xi$ onto the hemisphere $\SSS_\rad^\pm$ with the measure $d\eta'$ is measure preserving.

It follows (by Cavalieri's principle) that the bijective mapping $(v,\xi) \mapsto (v,\RRRR_v \xi)$ between the space $\Sigmapm$ with the measure $\mu$ (recall that it is defined by $d\mu(v,\xi) = |\langle v,\, n(\xi) \rangle|\, dv\, d\xi$) and the space $S^{d-1} \times \SSS_\rad^\pm$ with the measure defined by $dv\, d\eta'$ is also measure preserving. This implies, in particular, that
$$
\mu(\Sigmam) = \mu(\Sigmap) = |S^{d-1}|\, |B_\rad^{d-1}| = s_{d-1} b_{d-1} \rad^{d-1}.
$$

Now define the measures $\mu_\pm$ on $\SSS_\rad^\pm$ by $d\mu_\pm = s_{d-1} d\eta'$; then the mappings $\pi_\pm : \Sigmapm \to \SSS_\rad^\pm$ defined by
$$
\pi_\pm(v,\xi) = \RRRR_v \xi
$$
are measure preserving.

We have
\beq\label{eq5}
|\xi - \xi^+| = |\RRRR_v \xi - \RRRR_v \xi^+| \ge |\RRRR_v \xi - \RRRR_{v^+} \xi^+| - |\RRRR_v \xi^+ - \RRRR_{v^+} \xi^+|.
\eeq

From (\ref{Vge}) and (\ref{eq5}), taking into account that $|\RRRR_v \xi - \RRRR_{v^+} \xi^+| \le 2\rad$, we get
$$
V \ge \int_{\Sigmam\setminus\Sigma_\phio} |\RRRR_v \xi - \RRRR_v \xi^+|\, d\mu(v,\xi)
$$
$$
= \int_{\Sigmam\setminus\Sigma_\phio} |\RRRR_v \xi - \RRRR_{v^+} \xi^+|\, d\mu(v,\xi) - \int_{\Sigmam\setminus\Sigma_\phio} |\RRRR_v \xi^+ - \RRRR_{v^+} \xi^+|\, d\mu(v,\xi)
$$
\beq\label{Vge2}
\ge \int_{\Sigmam} |\RRRR_v \xi - \RRRR_{v^+} \xi^+|\, d\mu(v,\xi) - 2\rad\, \mu(\Sigma_\phio) - \int_{\Sigmam\setminus\Sigma_\phio} |\RRRR_v \xi^+ - \RRRR_{v^+} \xi^+|\, d\mu(v,\xi).
\eeq
Let us estimate the three terms in the right hand side of (\ref{Vge2}).

{\bf 3.1.} By Chebyshev's inequality one has $\mu(\Sigma_\phio) \le {\FFF(D)}/{f(\arccos\phio)}$, and therefore
\beq\label{eq1}
2\rad\, \mu(\Sigma_\phio) \le 2\rad\, \frac{\FFF(D)}{f(\arccos\phio)}.
\eeq

{\bf 3.2.} If $(v,\xi) \in \Sigmam\setminus\Sigma_\phio$, then the angle between $v$ and $v^+$ is less than $\phio$, and denoting by $\al = \al(v)$ and $\al^+ = \al(v^+)$ the angles formed by the vectors $v$ and $v^+$ with $(\bar 0, 1)$,\, $0 \le \al,\; \al^+ \le \pi$, we have $|\al - \al^+| \le \phio$.

If $d = 2$, one obviously has
$$
|\RRRR_v \xi^+ - \RRRR_{v^+} \xi^+| \le 2\rad\sin\frac{\phio}{2}.
$$
In the case $d \ge 3$ the estimate is more difficult.

Consider the 3-dimensional subspace spanned by the vectors $v$,\, $v^+$ and $(\bar 0, 1)$. The restrictions of $\RRRR_{v}$ and $\RRRR_{v^+}$ on this subspace are rotations by the angles $\al$ and $\al^+$, respectively. Let $w$ and $w^+$ be unit vectors in this subspace pointing at directions of the rotation axes. Both $w$ and $w^+$ are orthogonal to $(\bar 0, 1)$. The restriction of $\RRRR_{v^+}^{-1}\RRRR_{v}$ on this subspace acts as a rotation by an angle $\bt$, and its restriction on the orthogonal complement to this subspace is an identity. We have
$$
|\RRRR_v \xi^+ - \RRRR_{v^+} \xi^+| = |\RRRR_{v^+}^{-1}\RRRR_{v}\xi^+ - \xi^+| \le 2\rad\sin\frac{\bt}{2},
$$
therefore we need to estimate $\sin\frac{\bt}{2}$. To that end we shall proceed to some trigonometric calculations.

Introduce an orthonormal coordinate system  $x,\, y,\, z$ in the chosen subspace, where the third coordinate axis coincides with the $d$th axis of the original space $\RRR^d$ and the origin coincides with the origin in the space $\RRR^d$. In this system the coordinate vectors $v$,\, $v^+$, and $(0,\ldots,0,\, 1)$ take the form
$$
(\sin\al\cos\theta,\ \sin\al\sin\theta,\ \cos\al); \quad\ (\sin\al^+\cos\theta^+,\ \sin\al^+\sin\theta^+,\ \cos\al^+); \quad (0,\, 0,\, 1).
$$

One has
\beq\label{eqvv+}
\langle v, v^+ \rangle = \cos\al\, \cos\al^+ + \sin\al\, \sin\al^+ \cos(\theta-\theta^+).
\eeq
Further, one easily finds
that $w = (-\sin\theta,\ \cos\theta,\ 0)$,\, $w^+ = (-\sin\theta^+,\ \cos\theta^+,\ 0)$, and therefore
\beq\label{eqww+}
\langle w, w^+ \rangle = \cos(\theta-\theta^+).
\eeq
Taking into account that $\langle v, v^+ \rangle \ge \cos\phio$ and using (\ref{eqvv+}) and (\ref{eqww+}), one finds
\beq\label{wwge}
\langle w, w^+ \rangle \ge \frac{\cos\phio - \cos\al\, \cos\al^+}{\sin\al\, \sin\al^+}.
\eeq

In what follows we shall use the same notation $\RRRR_{v}$ and $\RRRR_{v^+}$ for the restrictions of the corresponding rotations on our 3D subspace. It is convenient to represent them in the quaternionic form: $\RRRR_{v}$ is the action $u \mapsto quq^{-1}$ of the quaternion
$$
q = \cos\frac{\al}{2} + \sin\frac{\al}{2}\ w,$$
and $\RRRR_{v^+}$ is the action $u \mapsto q_+ uq_+^{-1}$ of the quaternion
$$
q_+ = \cos\frac{\al^+}{2} + \sin\frac{\al^+}{2}\ w^+.
$$
Correspondingly, $\RRRR_{v^+}^{-1}\RRRR_{v}$ is the action of the quaternion
$$
q_+^{-1} q = \Big(\cos\frac{\al^+}{2} - \sin\frac{\al^+}{2}\ w^+\Big)\Big(\cos\frac{\al}{2} + \sin\frac{\al}{2}\ w\Big)
$$
$$
= \Big[ \cos\frac{\al}{2} \cos\frac{\al^+}{2} + \sin\frac{\al}{2} \sin\frac{\al^+}{2} \langle w, w^+ \rangle \Big] +
\Big[ \cos\frac{\al^+}{2}\sin\frac{\al}{2}\ w - \cos\frac{\al}{2}\sin\frac{\al^+}{2}\ w^+ + w^+ \times w \rangle \Big],
$$
which has the real part
$$
\cos\frac{\bt}{2} = \cos\frac{\al}{2} \cos\frac{\al^+}{2} + \sin\frac{\al}{2} \sin\frac{\al^+}{2} \langle w, w^+ \rangle.
$$
From this, using (\ref{wwge}) and utilizing  the double angle formulas for sine and cosine, one obtains the estimate
$$
\cos\frac{\bt}{2} \ge \cos\frac{\al}{2}\ \cos\frac{\al^+}{2} + \frac{\cos\phio - \cos\al\, (2\cos^2\frac{\al^+}{2} - 1)}{4\cos\frac{\al}{2}\ \cos\frac{\al^+}{2}}.
$$
Denoting $\cos\frac{\al^+}{2} =: z$, one comes to the inequality
$$
\cos\frac{\bt}{2} \ge \inf_{0\le z\le 1} \bigg( z \cos\frac{\al}{2} + \frac{\cos\phio - \cos\al\, (2z^2 - 1)}{4z\cos\frac{\al}{2}} \bigg).
$$
If $\al < \pi - \phio$, the infimum of the expression in the brackets is attained at $z = \sqrt{\cos\phio + \cos\al}/\sqrt 2$. Substituting this value in the latter inequality, one obtains
$$
\cos\frac{\bt}{2} \ge  \frac{\sqrt{\cos\phio + \cos\al}}{\sqrt 2\, \cos\frac{\al}{2}}.
$$
From here after some algebra one finally obtains
$$
\sin\frac{\bt}{2} \le  \frac{\sin\frac{\phio}{2}}{\cos\frac{\al}{2}}.
$$
Since we have $0 \le \al < \pi - \phio$, one easily concludes that the right hand side in this inequality is smaller than 1.

Thus, we have
$$
|\RRRR_v \xi^+ - \RRRR_{v^+} \xi^+| = |\RRRR_{v^+}^{-1} \RRRR_v \xi^+ -  \xi^+| \le 2\rad\, h_{\phio}(\al),
$$
where
$$
h_{\phio}(\al) = \left\{
\begin{array}{cl}
\frac{\sin\frac{\phio}{2}}{\cos\frac{\al}{2}}, & \text{if} \ 0 \le \al < \pi - \phio\\
1, & \text{if} \ \pi - \phio \le \al \le \pi
\end{array}
\right.
$$
in the case $d \ge 3$, and $h_{\phio}(\al) = \sin\frac{\phio}{2}$ in the case $d = 2$.
We now have an estimate for the last integral in the right hand side of (\ref{Vge2})
$$
\int_{\Sigmam\setminus\Sigma_\phio} |\RRRR_v \xi^+ - \RRRR_{v^+} \xi^+|\, d\mu(v,\xi) \le 2\rad \int_{\Sigmam} h_{\phio}(\al(v))\, d\mu(v,\xi).
$$

We now need to estimate the integral in the right hand side of this inequality. Integrating by $\xi$ gives us the factor $b_{d-1} \rad^{d-1}$. Integrating by $v$ over $S^{d-1}$ amounts to  integration with the differential $s_{d-2} \sin^{d-2}\al\, d\al$ over the interval $\al \in [0,\, \pi]$. Thus we get
$$
2\rad \int_{\Sigmam} h_{\phio}(\al(v))\, d\mu(v,\xi) = 2\rad^2 b_1 s_0 \int_0^\pi \sin\frac{\phio}{2}\, d\al = 8\rad^2 \pi\, \sin\frac{\phio}{2}
$$
in the case $d = 2$ and
$$
2\rad \int_{\Sigmam} h_{\phio}(\al(v))\, d\mu(v,\xi)
= 2\rad^d\, b_{d-1} s_{d-2} \Big( \sin\frac{\phio}{2} \int_0^{\pi-\phio} \frac{\sin^{d-2}\al}{\cos\frac{\al}{2}}\ d\al + \int_{\pi-\phio}^{\pi} \sin^{d-2}\al\, d\al \Big)
$$
in the case $d \ge 3$.

Introducing the functions $\IIII_d(\phi),\, \phi \in [0,\, \pi],\, d = 2, 3,\ldots$ by $\IIII_2(\phi) = \pi\sin\frac{\phi}{2}$ and
\beq\label{Id}
\IIII_d(\phi) = \sin\frac{\phi}{2} \int_0^{\pi-\phi} \frac{\sin^{d-2}\al}{\cos\frac{\al}{2}}\ d\al + \int_{\pi-\phi}^{\pi} \sin^{d-2}\al\, d\al.
\eeq
for $d \ge 3$, we can write
\beq\label{Vge3}
\int_{\Sigmam\setminus\Sigma_\phio} |\RRRR_v \xi^+ - \RRRR_{v^+} \xi^+|\, d\mu(v,\xi) \le 2\rad^d\, b_{d-1} s_{d-2}\, \IIII_d(\phio).
\eeq
For small values of $\phio$ we have the following asymptotic behavior:
\beq\label{Idapprox}
\IIII_2(\phio) = \frac{\pi}{2} \phio (1 + o(1)) \ \text{ and } \ \IIII_d(\phi) = 2^{d-3} B\big( \frac{d-1}{2}, \frac{d-2}{2} \big)\, \phi (1 + o(1)) \ (d \ge 3) \ \text{as} \ \phio \to 0^+.
\eeq
The function $\IIII_3$ can easily be calculated in the case $d=3$,\, $\IIII_3(\phio) = 4\sin\frac{\phio}{2} - 2\sin^2\frac{\phio}{2}$.


{\bf 3.3.} Now consider the first term in the right hand side of (\ref{Vge2}).

The mapping $T = T_{D,\rad}: \Sigmam \to \Sigmap$ preserves the measure $\mu$, and therefore induces a measure on $\Sigmam \times \Sigmap$ concentrated on the graph of $T$ and whose projections on $\Sigmam$ and $\Sigmap$ coincide with $\mu$. The push forward of this measure under the map\footnote{this map sends $(v,\xi,v^+,\xi^+)$ to $(\RRRR_v \xi,\, \RRRR_{v^+} \xi^+)$} $\pi_- \times \pi_+ : \Sigmam \times \Sigmap \to \SSS_\rad^- \times \SSS_\rad^+$ (let it be denoted by $\nu_D$) is a measure on $\SSS_\rad^- \times \SSS_\rad^+$ whose projections on $\SSS_\rad^-$ and on $\SSS_\rad^+$ coincide, respectively, with $\mu_-$ and $\mu_+$.

Therefore we have
$$
\int_{\Sigmam} |\RRRR_v \xi - \RRRR_{v^+}\xi^+|\, d\mu(v,\xi) = \int_{\SSS_\rad^- \times \SSS_\rad^+} |\eta - \eta^+|\, d\nu_D(\eta,\eta^+)
$$
$$
\ge \inf_\nu \int_{\SSS_\rad^- \times \SSS_\rad^+} |\eta - \eta^+|\, d\nu(\eta,\eta^+),
$$
where the infimum is sought among all measures $\nu$ whose projections on $\SSS_\rad^-$ and on $\SSS_\rad^+$ coincide with $\mu_-$ and $\mu_+$. That is, we have now a problem of mass transportation. This problem is actually easy to solve. We use the inequality $|\eta - \eta^+| \ge |\eta_d - \eta_d^+| = \eta_d^+ - \eta_d$ (since $\eta_d = - \sqrt{\rad^2 - \sum_{i=1}^{d-1}\eta_i^2} \le 0$ and $\eta_d^+ = \sqrt{\rad^2 - \sum_{i=1}^{d-1}(\eta_i^+)^2} \ge 0$) to get
$$
\int_{\SSS_\rad^- \times \SSS_\rad^+} |\eta - \eta^+|\, d\nu_D(\eta,\eta^+) \ge \int_{\SSS_\rad^- \times \SSS_\rad^+} \eta^+_d\, d\nu_D(\eta,\eta^+) - \int_{\SSS_\rad^- \times \SSS_\rad^+} \eta\, d\nu_D(\eta,\eta^+)
$$
$$
= \int_{\SSS_\rad^+} \eta^+_d\, d\mu_+(\eta^+) - \int_{\SSS_\rad^-} \eta_d\, d\mu_-(\eta)
= \int_{B_\rad^{d-1}} \Big( \rad^2 - \sum_{i=1}^{d-1}(\eta_i^+)^2 \Big)^{1/2}\, s_{d-1}\ d\eta_1^+ \ldots d\eta_{d-1}^+
$$
$$
-\int_{B_\rad^{d-1}} \bigg[ -\Big( \rad^2 - \sum_{i=1}^{d-1}\eta_i^2 \Big)^{1/2} \bigg] s_{d-1}\, d\eta_1 \ldots d\eta_{d-1} = s_{d-1} b_d\, \rad^d.
$$
(The value in the right hand side of this relation is really attained at the (optimal) measure supported on the subspace $(\eta_1,\ldots,\eta_{d-1}) = (\eta_1^+,,\ldots,\eta_{d-1}^+)$. This measure induces the transportation in the vertical direction sending any point of $\SSS_\rad^-$ to the point of $\SSS_\rad^+$ with the same abscissa.)

That is, we have
\beq\label{Vge4}
\int_{\Sigmam} |\RRRR_v \xi - \RRRR_{v^+}\xi^+|\, d\mu(v,\xi) \ge s_{d-1} b_d\, \rad^d.
\eeq

From (\ref{eq2}), (\ref{Vge2}), (\ref{eq1}), (\ref{Vge3}), and (\ref{Vge4}) we obtain
$$
s_{d-1} b_d\, \rad^d - s_{d-1} |D| \ge s_{d-1} b_d\, \rad^d - 2\rad\, \frac{\FFF(D)}{f(\phio)} - 2\rad^d\, b_{d-1} s_{d-2}\, \IIII_d(\phio).
$$
It follows that
\beq\label{eq6}
|D| \le \inf_{0 < \phio < \pi} \Big( \frac{2\rad\, \FFF(D)}{s_{d-1}}\, \frac{1}{f(\phio)} + \frac{2 b_{d-1} s_{d-2} \rad^d}{s_{d-1}}\, \IIII_d(\phio) \Big).
\eeq

Using the asymptotical formulas (\ref{f}) and (\ref{Idapprox}) for $f$ and $\IIII_d$, respectively, and replacing both terms in the right hand side of (\ref{eq6}) with their approximated values (as $\phio \to 0^+$), we obtain the expression
\beq\label{exp}
\frac{\al}{\kappao}\, \phio^{-k} + \bt \phio,
\eeq
where
$$
\al = \frac{2\kappao\rad\, \FFF(D)}{cs_{d-1}} \quad \text{and} \quad \bt =
\left\{
\begin{array}{ll}
2\rad^2, & \text{if}\, \ d = 2\\
\frac{1}{s_{d-1}} 2^{d-2} b_{d-1} s_{d-2} B\big( \frac{d-1}{2}, \frac{d-2}{2} \big) \rad^d, & \text{if}\, \ d \ge 3.
\end{array}
\right.
$$
The minimum of (\ref{exp}) is equal to $\frac{\kappao + 1}{\kappao}\, \al^{1/(\kappao+1)} \bt^{\kappao/(\kappao+1)}$ and is attained at $\phio_* = (\al/\bt)^{1/(\kappao+1)}$. Substituting this value $\phio_*$ in the right hand side of (\ref{eq6}) and raising both parts of the resulting inequality to the $(k+1)$th power, one obtains
\beq\label{eq7}
\left( \frac{|D|}{\rad^d} \right)^{\kappao+1} \le \frac{1}{c_d}\, \frac{\FFF(D)}{\rad^{d-1}}\, (1 + o(1)),
\eeq
where
$$
c_d = \left\{
\begin{array}{ll}
\frac{\kappao^\kappao}{(\kappao+1)^{\kappao+1}}\ \frac{\pi c}{2^\kappao}, & \text{if}\, \ d = 2\\
\frac{\kappao^\kappao}{(\kappao+1)^{\kappao+1}}\ \frac{cs_{d-1}^{\kappao+1} 2^{-1-d\kappao+2\kappao}}{\big( b_{d-1} s_{d-2} B\big( \frac{d-1}{2}, \frac{d-2}{2} \big) \big)^\kappao}, & \text{if}\, \ d \ge 3
\end{array}
\right.
$$
and $o(1)$ means a function of $\FFF(D)/\rad^{d-1}$ vanishing when its argument goes to zero. Reversing relation (\ref{eq7}), one gets
$$
\frac{\FFF(D)}{\rad^{d-1}} \ge c_d\, \left( \frac{|D|}{\rad^d} \right)^{\kappao+1}(1 + o(1));
$$
this time $o(1)$ means a function of $|D|/\rad^d$ vanishing when its argument goes to zero. Theorem \ref{t_main} is proved.

Let us now prove Theorem \ref{t2}. Here we have $f(\phio) = 1 - \cos\phio$ and use the notation $\fff$ in place of $\FFF$ in this particular case.

If $d = 2$, substitute $s_1 = 2\pi,\, b_1 = 2,\, s_0 = 2$, and $\IIII_2(\phio) = \pi\sin\frac{\phio}{2}$ into (\ref{eq6}) to obtain
$$
|D| \le \inf_{0 < \phio < \pi} \Big( \frac{\rad\fff(D)}{2\pi\sin^2 \frac{\phio}{2}} + 4\rad^2 \sin\frac{\phio}{2} \Big).
$$
Introducing the shorthand notation $z = \sin\frac{\phio}{2},\, A = |D|,\, \fff = \fff(D)$, this inequality can be rewritten as
\beq\label{inf}
A \le \inf_{0 < z < 1} \Big( \frac{\rad\fff}{2\pi z^2} + 4\rad^2 z \Big).
\eeq

Our goal is to prove the inequality
\beq\label{eqd=2}
A^3 \le \frac{54}{\pi}\, \rad^5 \fff,
\eeq
which is equivalent to statement (a) of Theorem \ref{t2}.

Consider two cases. If $\fff \le 4\pi\rad$, the infimum in (\ref{inf}) is attained at $z_* = (\fff/(4\pi\rad))^{1/3}$, and substituting $z_*$ in (\ref{inf}), we get (\ref{eqd=2}). On the other hand, if $\fff > 4\pi\rad$, we obviously have (since $D$ is contained in a circle of radius $\rad$)
$$
A^3 \le (\pi\rad^2)^3 < (6\rad^2)^3 = \frac{54}{\pi}\, \rad^5 \cdot 4\pi\rad < \frac{54}{\pi}\, R^5 \fff,
$$
and we again come to (\ref{eqd=2}). Thus, statement (a) of Theorem \ref{t2} is proved.

If $d = 3$, one has $s_2 = 4\pi,\, b_2 = \pi,\, s_1 = 2\pi$, and $\IIII_3(\phio) = 4\sin\frac{\phio}{2} - 2\sin^2\frac{\phio}{2}$, and inequality  (\ref{eq6}) takes the form
$$
|D| \le \inf_{0 < \phio < \pi} \Big[\frac{\rad\fff(D)}{4\pi\sin^2\frac{\phio}{2}} + 2\pi\rad^3 \Big( 2\sin\frac{\phio}{2} - \sin^2\frac{\phio}{2} \Big)\Big].
$$
Introducing the notation $\tilde A = |D|/(2\pi\rad^3), \ \tilde\fff = \fff(D)/(8\pi^2\rad^2)$, and $z = \sin\frac{\phio}{2}$, one rewrites the last inequality in the form
\beq\label{tildeA}
\tilde A \le \inf_{0 \le z \le 1} h(z), \quad \text{where} \ \, h(z) = \frac{\tilde\fff}{z^2} + 2z - {z^2}.
\eeq

We are going to prove the inequality
\beq\label{ind=3}
\tilde A^3 \le 27 \tilde\fff,
\eeq
which is equivalent to statement (b) of Theorem \ref{t2}.

After a simple algebra one concludes that if $\tilde \fff > 27/256$, we have $h'(z) < 0$ for all $z > 0$.    has a unique zero $z = 3/4$. If $0 < \tilde \fff \le 27/256$, the equation $h'(z) = 0$ has two positive zeros (coinciding when $\tilde\fff = 27/256$). The smallest zero $z_* = z_*(\tilde\fff)$ (which is a local minimizer of $h$ if $\fff$ is strictly smaller than $27/256$) satisfies the inequality $0 < z_* \le 3/4$. It is also straightforward to check that
\beq\label{tildeF}
\tilde\fff = z_*^3 (1 - z_*).
\eeq

Consider two cases. If $\tilde\fff \le 27/256$, we substitute $z = z_*$ in (\ref{tildeA}) and use (\ref{tildeF}) to obtain $\tilde A \le 3z_* (1 - 2z_*/3)$. Taking the third power of both sides of this inequality and using that $(1 - 2z/3)^3 < 1 - z$ for $0 < z \le 3/4$, we come to (\ref{ind=3}). If, otherwise, $\tilde\fff > 27/256$, we use that $|D| \le \frac 43 \pi\rad^3$, and therefore $\tilde A \le 2/3$. It follows that $\tilde A^3 \le 8/27 < 27 \cdot 27/256 < 27 \tilde\fff$, and (\ref{ind=3}) again follows. Thus, statement (b) of Theorem \ref{t2} is also proved.

\end{document}